\documentclass{article}
\usepackage{amsfonts}
\usepackage{amsthm}
\usepackage{amsmath}
\usepackage{amscd}
\usepackage{authblk}
\title{Irrotational Stokes flow in three system of coordinates.}
\author{Eleftherios Protopapas}
\affil{National Technical University of Athens, School of Applied Mathematical and Physical Sciences, Department of Mathematics.}
\affil{lprotopapas@math.ntua.gr}
\date{}

\begin{document}

\maketitle

\begin{abstract}
Irrotational flow is described with the second order elliptic partial differential equation $E^2\psi=0,$ where $\psi$ is the  function to be derived and $E^2$ is the Stokes operator. In the present paper we derive the solution of $E^2\psi=0$ in three axisymmetric system of coordinates: the parabolic, the tangent sphere and the cardioid. We prove that Stokes equation separates variables in the parabolic coordinate system, but it R-separates variables  in the other two. At the end of this manuscript we summarize all the known results for the 0-eigenspace of $E^2$ in axisymmetric system of coordinates.\\
\textbf{Keywords:} Irrotational Stokes flow, 0-eigenspace, Stokes operator, Parabolic, Tangent Sphere, Cardioid.
\end{abstract}



\section{Introduction.}
Stokes operator, $E^2,$ is a second order partial elliptic operator and the equation $E^2\psi=0,$ describes the creeping axisymmetric irrotational flow  \cite{happelbrenner}.  Steady creeping flow is the flow where the inertial forces dominate over the viscous ones and it is described through the forth order partial elliptic equation $E^4\psi=0,$ where $\psi$ is the stream function and $E^4=E^2oE^2$ is the Stokes bistream operator \cite{happelbrenner}. This is used to model flow phenomena in chemical engineering, biology, medicine etc. in cases of low Reynolds number. Indicatively, in previous studies we used Stokes flow to model blood flow. Dassios et al. \cite{blood} develop an analytical model for the blood's plasma flow past a red blood cell, while Hadjinicolaou et al. \cite{sedi} expand this model to describe the sedimentation of a red blood cell. Hadjinicolaou \cite{ldl1} and Hadjinicolaou, Protopapas \cite{ldl2} studied the relative motion of two aggregated low density lipoproteins in blood's plasma. Recently  Hadjinicolaou and Protopapas \cite{hpkuwa} modeled blood's plasma flow through a swarm of red blood cells as a Stokes flow problem, providing analytical solutions. Extensive studies of Stokes flow in chemical engineering are also known. Pitter et al. \cite{pitpruham} obtained accurate solutions for the flow past a thin oblate spheroid in order to fill the need for collision efficiencies of ice crystals. The flow in a fluid cell contained between two confocal spheroidal surfaces with Kuwabara-type boundary conditions deriverd by Dassios et al. in \cite{semi}. Deo \cite{deo1} studied Stokes flow past a swarm of porous circular cylinders with Happel and Kuwabara boundary conditions, while the flow through a swarm of porous nanocylindrical particles enclosing a solid cylindrical core using Kuwabara type boundary conditions was discussed by Deo and Yadav \cite{deoyadav}. \\
The closed form solution of $E^4\psi=0$ is known in many axisymmetric system of coordinates given as series expansion of  eigenfunctions. In the spherical system the equation separates variables \cite{happelbrenner}, while in the spheroidal system of coordinates it semi-separates variables \cite{semi}. Moreover in the inverted spheroidal systems it R-semiseparates variables \cite{rsemi1}, \cite{rsemi2}. In each one of these coordinate systems the first step for obtaining the solution of $E^4\psi=0$ is the derivation of the 0-eigenspace of $E^2$ and then its generalized 0-eigenspace. This way the function that expresses the irrotational flow is then consisted of these eigenfunctions. Furthermore the creeping axisymmetric irrotational flow has been also solved in bispherical and in toroidal coordinate systems \cite{deotiwari}. 
\\ In the present study we derive the solution of the irrotational flow in the parabolic, the tangent sphere and the cardioid coordinate systems. These systems must be employed in every case that one of the coordinate surface matches the shape of the particles that are employed in the flow. We use the method of separation of variables for the parabolic system, in which 
the 0-eigenspace employ Bessel and modified Bessel functions of the first order. In the tangent sphere coordinate system equation $E^2\psi=0$ R-separates variables, where R is the inverse of the euclidean distance and the employed functions are the sine, cosine and modified Bessel of the first order. Moreover in cardioid geometry Stokes equation also R-separates variables, with R being the inverse of the euclidean distance and the involved functions are the Bessel and modified Bessel of the first order.\\
The structure of this manuscript is as follows: in section 2 we present the mathematical background about the rotational systems and the Stokes operator, in sections 3, 4, 5 we derive the 0-eigenspaces of Stokes operator in the aforementioned coordinate systems and in section 6 we discuss our results.

\section{Mathematical background.}
\label{back}
A rotational system of coordinates  $(q_1,q_2,\varphi), \; \varphi \in [0, 2\pi)$ is defined \cite{happelbrenner} as
\begin{equation}
\label{eqn: coord}
\left\{ \begin{array}{*{20}{c}}
{x=\rho(q_1,q_2) cos \varphi}\\
{y=\rho(q_1,q_2) sin \varphi}\\
{z=z(q_1,q_2)}
\end{array}  \right., 
\end{equation}
where the metric coefficients are
\begin{equation}
\label{eqn: h1}
h_1=\frac{1}{\sqrt{\left(\displaystyle{ \frac{\partial \rho}{\partial q_1}}\right)^2+\left(\displaystyle{ \frac{\partial z}{\partial q_1}}\right)^2}},
\end{equation}
\begin{equation}
\label{eqn: h2}
h_2=\frac{1}{\sqrt{\left(\displaystyle{ \frac{\partial \rho}{\partial q_2}}\right)^2+\left(\displaystyle{ \frac{\partial z}{\partial q_2}}\right)^2}}
\end{equation}
and the radial cylindrical coordinate is 
\begin{equation}
\label{eqn: omega}
\varpi=|\rho(q_1,q_2)|.
\end{equation}
Stokes operator in the axisymmetric system of coordinates $(q_1,q_2,\varphi)$ has the form
\begin{equation}
\label{eqn: stokesop}
 E^2=h_1h_2 \varpi \left[\frac{\partial}{\partial q_1}\left(\frac{h_1}{h_2\varpi} \frac{\partial}{\partial q_1}\right)+\frac{\partial}{\partial q_2}\left(\frac{h_2}{h_1\varpi} \frac{\partial}{\partial q_2} \right) \right]
\end{equation}
and if $\psi=\psi(q_1,q_2)$ is  function that describes the irrotational flow, Stokes equation $E^2\psi=0$ is written as
\begin{equation}
\label{eqn: stokeseq}
 \frac{\partial }{\partial q_1} \left(\frac{h_1}{h_2\varpi} \frac{\partial \psi}{\partial q_1}\right)+\frac{\partial }{\partial q_2}\left(\frac{h_2}{h_1\varpi} \frac{\partial\psi}{\partial q_2} \right)=0.
\end{equation}

\section{Eigenfunctions of Stokes operator in the parabolic system of coordinates.}
In parabolic coordinate system \cite{moonspencer} every point $(x,y,z)$ in the Cartesian coordinate system is expressed with $(\mu ,\nu,\varphi ),$ where $\mu,\nu \geq 0$ and
\begin{equation}
\label{eqn: parabolic} 
\left\{ {\begin{array}{*{20}{c}}
{x = \mu \nu \cos (\varphi )}\\
{y = \mu \nu \sin (\varphi )}\\
{z = \displaystyle{\frac{{{\mu ^2} - {\nu ^2}}}{2}}}
\end{array}} \right.,
\end{equation}
while Stokes operator assumes the form
\begin{equation}
\label{eqn: E2parabolic}
E^2=\frac{1}{\mu^2+\nu^2}\left(\frac{\partial^2}{\partial \mu^2}-\frac{1}{\mu}\frac{\partial}{\partial \mu} -\frac{1}{\nu}\frac{\partial}{\partial \nu} +\frac{\partial^2}{\partial \nu^2}\right).
\end{equation}
In equation $E^2\psi_p(\mu,\nu)=0,$ we set \begin{equation}
\label{eqn: paraxwri}
\psi_p(\mu,\nu)=M(\mu)N(\nu),
 \end{equation}
deriving that
\begin{equation}
\label{eqn: paraxwri}
\frac{M''}{M}-\frac{1}{\mu}\frac{M'}{M}=-\frac{N''}{N}+\frac{1}{\nu}\frac{N'}{N},
\end{equation}
which denotes that $E^2\psi=0$ separates variables. If $\lambda=n^2>0$ is the separation constant from \eqref{eqn: paraxwri} yields 
\begin{equation}
\label{eqn: paraxwriM}
\mu^2 M''-\mu M'-n^2\mu^2M=0
\end{equation}
and
\begin{equation}
\label{eqn: paraxwriN}
\nu^2 N''-\nu N'+n^2\nu^2N=0,
\end{equation}
which are the two ODEs that we need to solve in order to derive the 0-eigenspace of Stokes operator in the parabolic system of coordinates.
In Appendix A we derive analytically the solution of \eqref{eqn: paraxwriM} and in Appendix B the solution of \eqref{eqn: paraxwriN}, which are
\begin{equation}
\label{eqn: Mxwripara}
M(\mu)=c_1\mu I_1(n \mu)+c_2 \mu K_1(n \mu),
\end{equation}
\begin{equation}
\label{eqn: Nxwripara}
N(\nu)=c_3\nu J_1(n \nu)+c_4 \nu Y_1(n \nu),
\end{equation}
where $J_1,Y_1$ are Bessel functions of the first order and of the first and the second kind respectively, $I_1,K_1$ are modified Bessel functions of the first order and of both kinds \cite{lebe}, while $c_1,c_2,c_3,c_4$ are constants. Therefore
\begin{equation}
\label{psipara}
\psi_p(m,n)=\sum_{n=1}^{\infty} \left[A_n\mu I_1(n \mu)+B_n \mu K_1(n \mu)\right]\left[C_n\nu J_1(n \nu)+D_n \nu Y_1(n \nu)\right],
\end{equation}
where $A_n,B_n,C_n,D_n$ are constants to be calculated  from the boundary conditions.
\section{Eigenfunctions of Stokes operator in the tangent sphere system of coordinates.}
In tangent sphere coordinate system \cite{moonspencer} every point $(x,y,z)$ in the Cartesian coordinate system is expressed with $(\mu ,\nu,\varphi ),$ where $\mu > 0,\; \nu \in \mathbb{R}$ and 
\begin{equation}
\label{eqn: tangent} 
\left\{ {\begin{array}{*{20}{c}}
{x = \displaystyle{\frac{\mu \cos(\varphi ) }{{{\mu ^2} + {\nu ^2}}}}}  \vspace{1 mm} \\
{y = \displaystyle{\frac{\mu \sin (\varphi ) }{{{\mu ^2} + {\nu ^2}}}}}  \\ \vspace{1 mm}
{z = \displaystyle{\frac{\nu }{{{\mu ^2} + {\nu ^2}}}}}\\
\end{array}}  \right.,
\end{equation}
while Stokes operator assumes the form
\begin{equation}
\label{eqn: E2tangent}
E^2=(\mu^2+\nu^2)^2\left[\frac{\partial^2}{\partial \mu^2}+\frac{\mu^2-\nu^2}{\mu(\mu^2+\nu^2)}\frac{\partial}{\partial \mu} +\frac{2\nu}{\mu^2+\nu^2} \frac{\partial}{\partial \nu} +\frac{\partial^2}{\partial \nu^2}\right].
\end{equation}
In equation $E^2\psi_t(\mu,\nu)=0,$ we set 
\begin{equation}
\label{eqn: Rpsitang}
\psi_t(\mu,\nu)=\frac{1}{\sqrt{\mu^2 + \nu^2}}M(\mu)N(\nu),
\end{equation}
deriving that
\begin{equation}
\label{eqn: tangxwri}
\frac{M''}{M}-\frac{1}{\mu}\frac{M'}{M}=-\frac{N''}{N}=\lambda,
\end{equation}
which is in separable form and $\lambda$ is the separation constant. We assume that $\lambda=n^2>0$ and we solve the O.D.Es. that arise from \eqref{eqn: tangxwri}, deriving that
\begin{equation}
\label{eqn: Mxwritang}
M(\mu)=c_1\mu I_1(n \mu)+c_2 \mu K_1(n \mu)
\end{equation}
and
\begin{equation}
\label{eqn: Nxwritang}
N(\nu)=c_3\cos(n \nu)+c_4 \sin(n \nu).
\end{equation}
From these results it yields that $E^2\psi_t(\mu,\nu)=0,$ R-separates variables, where 
\begin{equation}
\label{eqn: Rpsitang}
R(\mu,\nu)=\sqrt{\mu^2 + \nu^2},
\end{equation}
which is the inverse of the euclidean distance and 
\begin{equation}
\label{psitang}
\psi_t(m,n)=\frac{1}{\sqrt{\mu^2 + \nu^2}}\sum_{n=1}^{\infty} \left[A_n\mu I_1(n \mu)+B_n \mu K_1(n \mu)\right]\left[C_n \cos(n \nu)+D_n \sin(n \nu)\right],
\end{equation}
where $A_n,B_n,C_n,D_n$ are constants to be calculated from the boundary conditions.

\section{Eigenfunctions of Stokes operator in the cardioid system of coordinates.}
In cardioid coordinate system \cite{moonspencer} every point $(x,y,z)$ in the Cartesian coordinate system is expressed with $(\mu ,\nu,\varphi ),$ where $\mu,\nu \geq 0$ and
\begin{equation}
\label{eqn: cardioid} 
\left\{ {\begin{array}{*{20}{c}}
{x = \displaystyle{\frac{{\mu \nu \cos(\varphi ) }}{{{{\left( {{\mu ^2} + {\nu ^2}} \right)}^2}}} }}\\ \vspace{1 mm}
{y = \displaystyle{\frac{{\mu \nu \sin (\varphi )}}{{{{\left( {{\mu ^2} + {\nu ^2}} \right)}^2}}} }}\\ \vspace{1 mm}
{z = \displaystyle{\frac{{{\mu ^2} - {\nu ^2}}}{{2{{\left( {{\mu ^2} + {\nu ^2}} \right)}^2}}}}}
\end{array}} \right.,
 \end{equation}
while Stokes operator assumes the form
\begin{equation}
\label{eqn: E2cardio}
E^2=(\mu^2+\nu^2)^3\left[\frac{\partial^2}{\partial \mu^2}+\frac{3\mu^2-\nu^2}{\mu(\mu^2+\nu^2)}\frac{\partial}{\partial \mu} +\frac{3\nu^2-\mu^2}{\nu(\mu^2+\nu^2)}\frac{\partial}{\partial \nu} +\frac{\partial^2}{\partial \nu^2}\right].
\end{equation}
In equation $E^2\psi_c(\mu,\nu)=0,$ we set
\begin{equation}
\label{eqn: Rsepcardio}
\psi_c(\mu,\nu)=\frac{1}{\sqrt{2}\left(\mu^2+\nu^2\right)}M(\mu)N(\nu),
\end{equation}
deriving that
\begin{equation}
\label{eqn: cardioxwri}
\frac{M''}{M}-\frac{1}{\mu}\frac{M'}{M}=-\frac{N''}{N}+\frac{1}{\nu}\frac{N'}{N}=\lambda,
\end{equation}
which is in separable form and $\lambda$ is the separation constant.
\\Solving the O.D.Es. that arise from \eqref{eqn: cardioxwri} with $\lambda=n^2>0$ we have
\begin{equation}
\label{eqn: Mxwricardio}
M(\mu)=c_1\mu I_1(n \mu)+c_2 \mu K_1(n \mu)
\end{equation}
and
\begin{equation}
\label{eqn: Nxwrixcardio}
N(\nu)=c_3\nu J_1(n \nu)+c_4 \nu Y_1(n \nu).
\end{equation}
From these solutions it yields that $E^2\psi_c(\mu,\nu)=0,$ R-separates variables, where 
\begin{equation}
\label{eqn: Rpsicardio}
R(\mu,\nu)=\sqrt{2} \left(\mu^2 + \nu^2 \right),
\end{equation}
which is the inverse of the euclidean distance and 
\begin{equation}
\label{psicardio}
\psi_c(m,n)=\frac{1}{\mu^2 + \nu^2}\sum_{n=1}^{\infty} \left[A_n\mu I_1(n \mu)+B_n \mu K_1(n \mu)\right]\left[C_n\nu J_1(n \nu)+D_n \nu Y_1(n \nu)\right],
\end{equation}
where $A_n,B_n,C_n,D_n$ are constants to be calculated from the boundary conditions.

\section{Discussion.}

In the present paper we derived the exact solutions of the irrotational axisymmetric creeping flow in the parabolic, the tangent-sphere and the cardioid system of coordinates. The solution in the parabolic system of coordinates is obtained in a simply separable from and consists of Bessel functions and modified Bessel functions of the first order. In tangent sphere geometry, the equation $E^2\psi=0$ R-separates variables, with R being the inverse of the euclidean distance and the 0-eigenspace consists of sine, cosine functions and modified Bessel functions of the first order. Moreover in cardioid coordinate system Stokes equation R-separates variables with R being the inverse of the euclidean distance and the 0-eigenspace consists of Bessel functions and modified Bessel functions of the first order.
\\
\appendix
\section{Solution of $x^2y''(x)-xy'(x)-n^2x^2y(x)=0, \;x \neq 0.$} 
If we set $$y(x)=xu(nx),$$ it yields $$y'(x)=u(nx)+nxu'(nx)$$ and $$y''(x)=2nu'(nx)+n^2xu''(nx),$$ the equation becomes
$$n^2x^2u''(nx)+nxu'(nx)-(n^2x^2+1)u(nx)=0.$$
Substituting $t=nx,$ we get
$$t^2u''(t)+tu'(t)-(t^2+1)u(t)=0,$$
which is a modified Bessel equation of the first order \cite{lebe}, therefore the general solution is
$$u(t)=c_1I_1(t)+c_2K_1(t),$$
so 
$$y(x)=c_1xI_1(nx)+c_2xK_1(nx),$$
where $I_1,K_1$ are modified Bessel functions of the first and the second kind respectively.

\section{Solution of $x^2y''(x)-xy'(x)+n^2x^2y(x)=0,\;x \neq 0.$} 
If we set $$y(x)=xu(nx),$$ it yields $$y'(x)=u(nx)+nxu'(nx)$$ and $$y''(x)=2nu'(nx)+n^2xu''(nx),$$ the equation becomes
$$n^2x^2u''(nx)+nxu'(nx)+(n^2x^2-1)u(nx)=0.$$
Substituting $t=nx,$ we get
$$t^2u''(t)+tu'(t)+(t^2-1)u(t)=0,$$
which is a Bessel equation of the first order \cite{lebe}, therefore the general solution is
$$u(t)=c_1J_1(t)+c_2Y_1(t),$$
so 
$$y(x)=c_1xJ_1(nx)+c_2xY_1(nx),$$
where $J_1,Y_1$ are Bessel functions of the first and the second kind respectively.

\end{document}